\documentclass[11pt]{article}

\usepackage{microtype}
\usepackage{graphicx}
\usepackage{subfigure}
\usepackage{booktabs} 
\usepackage{amssymb}
\usepackage{amsmath}
\usepackage{multirow}
\usepackage{mathtools}
\usepackage[noend]{algorithmic}
\usepackage{algorithm}

\usepackage{siunitx}
\usepackage{color}
\usepackage{ulem}
\usepackage{sidecap}
\usepackage{wrapfig}
\usepackage{lipsum}

\DeclareMathOperator*{\argmin}{arg\,min}

\newcommand{\defeq}{\vcentcolon=}

\usepackage{microtype} 
\usepackage{booktabs}  

\usepackage{amsmath}
\usepackage{amsthm}

\usepackage[final]{automl_arxiv}
\usepackage{automl_arxiv}

\usepackage{natbib}

\usepackage{filecontents}

\title{Efficient Meta Subspace Optimization}

\author[1]{\nameemail{Yoni Choukroun}{choukroun.yoni@gmail.com}}
\author[1]{\nameemail{Michael Katz}{michael.c.katz@gmail.com}}

\affil[1]{Huawei Technologies \\ Israel Research Center}

\hypersetup{%
  pdfauthor={Y Choukroun, M Katz}, 
  pdftitle={Efficient Meta Subspace Optimization},
  pdfsubject={Efficient Meta Subspace Optimization}
}

\begin{document}

\abovedisplayskip=12pt plus 3pt minus 9pt
\belowdisplayskip=12pt plus 3pt minus 9pt
\abovedisplayshortskip=0pt plus 3pt
\belowdisplayshortskip=7pt plus 3pt minus 4pt
\abovedisplayskip=6pt
\belowdisplayskip=6pt
\abovedisplayshortskip=-6pt
\belowdisplayshortskip=1pt

\maketitle

\begin{abstract}
Subspace optimization methods have the attractive property of reducing large-scale optimization problems to a sequence of low-dimensional subspace optimization problems.
However, existing subspace optimization frameworks adopt a fixed update policy of the subspace and therefore appear to be sub-optimal.
In this paper, we propose a new \emph{Meta Subspace Optimization} (MSO) framework for large-scale optimization problems, which allows to determine the subspace matrix at each optimization iteration. 
In order to remain invariant to the optimization problem's dimension, we design an \emph{efficient} meta optimizer based on very low-dimensional subspace optimization coefficients, inducing a rule-based method that can significantly improve performance. 
Finally, we design and analyze a reinforcement learning (RL) procedure based on the subspace optimization dynamics whose learnt policies outperform existing subspace optimization methods.
\end{abstract}

\section{Introduction}
\label{introduction}
First order optimization methods are ubiquitous in numerical optimization. Gradient based frameworks are computationally efficient and are proven to converge in various
settings.
However, these methods suffer from inherent slowness due to the poor local linear approximation which is extremely inefficient in ill-conditioned problems.

One popular approach to address these problems is to resort to second-order methods \citep{NocedWrite}, which can converge extremely fast to stationary points, especially in the proximity of the solution where the problem has a good quadratic approximation.
However, when the number of variables is very large, there is a need for optimization algorithms whose storage requirement and computational cost per iteration grow at most linearly with the problem's dimensions. This constraint led to the development of a broad family of methods e.g., variable metric methods and \emph{subspace optimization}. 

Subspace optimization iteratively searches for the optimum in a low-dimensional subspace spanned by vectors obtained from a first order oracle \citep{conn1996iterated}.
The vectors are generally updated at each iteration by \emph{axiomatically} removing the oldest vector of the subspace obtained during the iterative optimization.

Meta optimization is an emerging field of research, focused on \emph{learning} optimization algorithms by training a parametrized optimizer on a distribution of tasks \citep{andrychowicz2016learning}. 
These meta optimizers (MOs) extend the handcrafted update rules of classical first order methods with high-dimensional parametrized models which provide better updates of the optimizee and have been shown to outperform classical baselines in restricted settings.

The subspace optimization paradigm has been proven to be a method of choice for deterministic large scale problems which is the focus of this paper.
In this work, we extend the idea of learned optimization to the paradigm of subspace optimization.
{\color{black}Unlike existing techniques which axiomatically remove the oldest direction, the proposed} meta-optimizer is trained to predict the best direction to be removed from the subspace given a low-dimensional state obtained from previous subspace optimization steps.  
The main contributions of the paper are as follows.
\begin{itemize}
    \item We extend for the first time the idea of meta optimization to the field of subspace optimization.
    \item We propose a meta-optimization framework which includes processing of very low-dimensional information obtained from the subspace spanning coefficients, enabling highly efficient training and deployment of the meta optimizer.
    \item We design a rule-based approach that can outperform popular subspace optimization methods, and 
     propose a RL agent able to determine the optimal subspace matrix and obtain SOTA results.
     \item Finally, we analyse the agent's actions throughout the optimization, and demonstrate its ability to learn optimal policies \emph{adapted} to the given tasks.
\end{itemize}

\section{Background}
\label{related_works}
\subsection{Subspace Optimization}
\label{subsec:Background-subspace-opt}
The core idea of subspace optimization is to perform the optimization of the objective function in a small subspace spanned by a set of directions obtained from an available oracle.
Denoting a function $f:\mathbb{R}^{n} \to \mathbb{R}, f\in C^{2}$ to be minimized and $P_k\in \mathbb{R}^{n\times d}, d \ll n$ as the set of $d$ directions at iteration $k$, an iterated subspace optimization method aims at solving the following inner minimization problem
\begin{equation}\label{eq:sesop_minimization}
\begin{aligned}
	\alpha_k = \argmin_{\alpha \in \mathbb{R}^d} f(x_k + P_k \alpha),
\end{aligned}
\end{equation}
followed by the update rule
\begin{equation}
    x_{k+1} = x_k +P_k \alpha_k \label{eq: subspace optimization update rule}.
\end{equation}
The dimensions of the problem are then reduced from the entire optimization space  $\mathbb{R}^n$, to a controlled $d$-dimensional subspace of $\mathbb{R}^n$ spanned by the columns of $P_k$. 

A direct benefit of subspace optimization is that the low-dimensional optimization task at every iteration can be addressed efficiently using heavier optimization tools such as second order methods.
The main computational burden here is the need to multiply the spanning directions by second order derivatives, which can be implemented efficiently using Hessian-vector product rules, i.e. $(\partial^{2}f)\cdot v=\partial(\partial f \cdot v)$.
Alternatively, quasi-Newton methods can be used for the inner optimization.

The subspace structure may vary depending on the chosen optimization technique.
Early methods proposed to extend the minimization to a $d$-dimensional subspace spanned by $d$ various previous directions, such as gradients, conjugate directions, previous outer iterations or Newton directions \citep{cragg1969study,miele1969study,dennis1987generalized,conn1996iterated}.
The Krylov descent method defines the subspace as $\text{span}\{H^{0}\nabla f,...,H^{d-1}\nabla f\}$ for some preconditioning matrix $H$, for example $H=\nabla^2f$ in \citep{vinyals2012krylov}.

Related to Krylov subspaces, the Conjugate Gradient (CG) method \citep{hestenes1952methods} reduces the search space to current gradient and previous step, i.e. $\text{span}\{p_{k},\nabla f(x_{k})\}$, where $p_{k}=x_{k}-x_{k-1}$. CG possesses a remarkable linear convergence rate compared to steepest descent methods in the quadratic case, related to the expanding manifold property \citep{NocedWrite}.
\cite{nemirovski1982orth} provided optimal worst case complexity in the convex setting with the ORTH-method by defining the subspace as $\text{span}\{\sum_{j=0}^{k}w_{j}\nabla f(x_{j}),x_{k}-x_{0},\nabla f(x)\}$, with appropriate weights $\{w_{j}\}_{j=0}^{k}$.

The Sequential Subspace Optimization (SESOP) algorithm \citep{Narkiss-2005} extends ORTH by adding the previous search directions $\{p_{k-i}\}_{i=0}^{d-3}$. This way, the method generalizes the CG method by allowing truncated approximation of the expanding manifold property on non-linear objectives, coupled with the worst case optimality safeguard of the ORTH method.
Of course any other (valuable) direction can be further embedded into the subspace structure in order to improve convergence as in \citep{conn1996iterated,ZibEladSPM,seboost,choukroun2020primal}. 

\subsection{Meta Optimization}
While deep learning has achieved tremendous success by replacing hand-crafted feature engineering with automatic feature learning from large amounts of data, hand-crafted optimization algorithms, such as momentum and adaptive gradient based methods \citep{sutskever2013importance,kingma2014adam} are still mainly in use for training neural networks.
These optimizers typically require careful tuning of hyper parameters
and extensive expert supervision in order to be used effectively for different model architectures and datasets \citep{choi2019empirical}. 
Meta optimization attempts to replace hand-designed optimizers by training a parametrized optimizer on a set of tasks and then applying it to the optimization of different tasks. 

One meta optimization approach aims to train a controller to automatically adapt the hyper parameters of an existing hand-crafted optimizer based on the training dynamics. 
\citep{daniel2016learning} exploits features based on the variance of the predicted and observed changes in function value with a linear policy mapping of the learning rate for SGD and RMSProp optimizers. 
\citep{xu2017reinforcement} learns a small LSTM \citep{hochreiter1997long} network using the current training loss as input in order to predict the learning rate. 
\citep{xu2019learning} scales the learning rate using an LSTM or an MLP based on features extracted from the training dynamics such as the train and validations losses, the variance of the network predictions, and the last layer's parameters statistics.




\begin{wrapfigure}{l}{0.47\textwidth}
\vskip -15pt
\begin{center}
\includegraphics[trim={0 50 0 150}, width=0.47\textwidth]{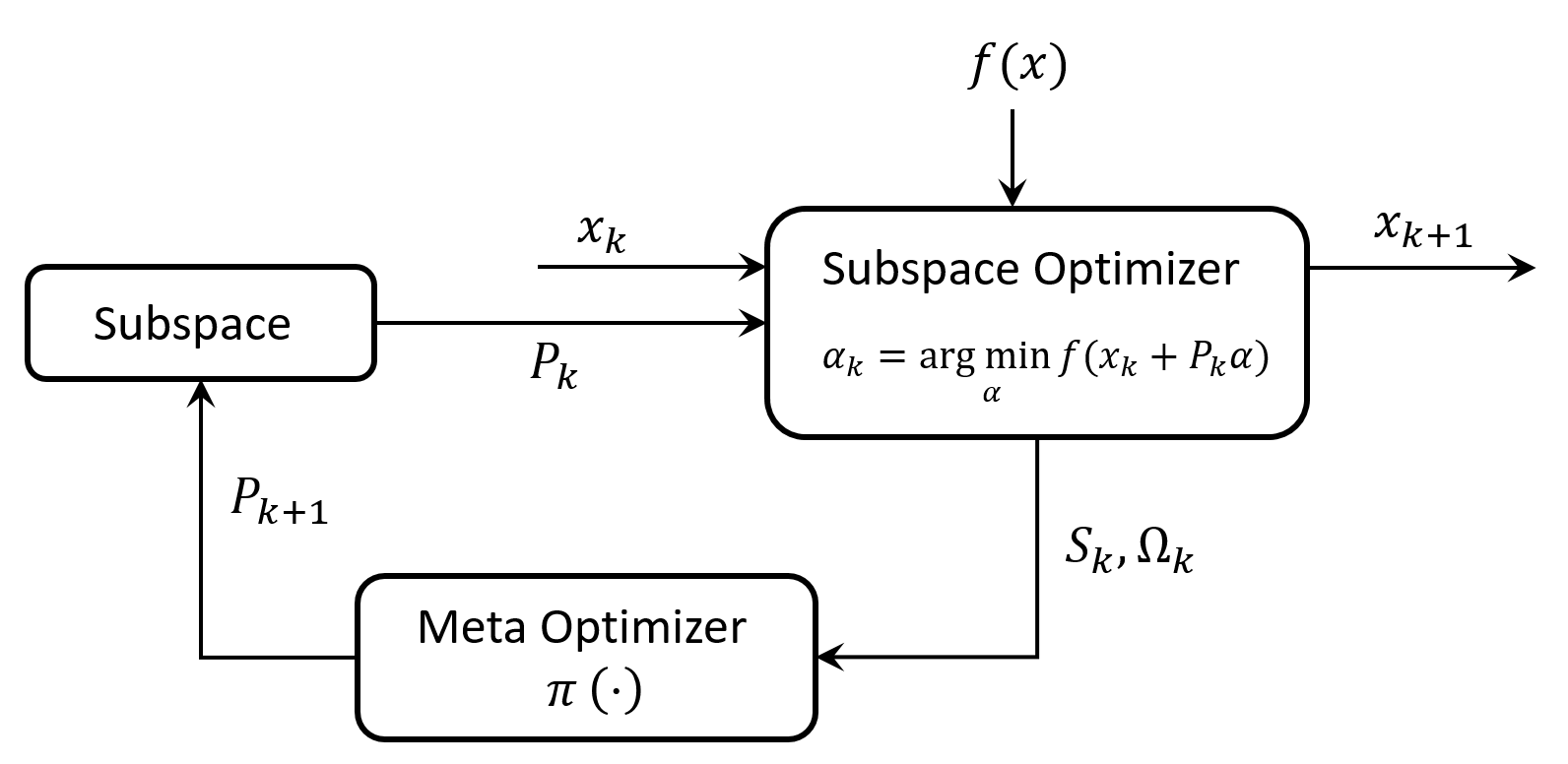}
\end{center}
\caption{Diagram of our proposed framework.
The subspace optimizer processes the current iterate $x_{k}$ and the subspace matrix $P_{k}$ and determines the next iterate: the set of directions $S_{k}$ and the subspace optimization information $\Omega_{k}$. The MO $\pi$ determines the next iteration's subspace matrix.}
\vskip -10pt
\label{fig:diagram}
\end{wrapfigure}

Another thread of research is devoted to the development of more expressive learned optimizers meant to replace existing optimizers entirely. 
\cite{andrychowicz2016learning} propose a learned update rule given by \mbox{$\theta_{t+1} = \theta_t + g_t\left(\nabla f\left(\theta\right), \phi\right)$} where $g_t$ is a \emph{coordinate-wise} LSTM that maps the optimizee's gradients to the new step vector.
\cite{li2017learning} proposed a similar approach which relies on policy search to compute the meta-parameters of the optimizer. 
Subsequent studies have improved the robustness of these approaches.
\cite{wichrowska2017learned} introduces a hierarchical RNN architecture which captures inter-parameter dependencies and uses an ensemble of small tasks with diverse loss landscapes. 
\cite{lv2017learning} incorporates random scaling of the optimizee when training the RNN-based optimizer which is shown to improve its generalization. 
\cite{metz2019understanding} improve the meta-training of the learned optimizer by constructing two different unbiased gradient estimators for the variational loss on optimizer performance. 
Finally, \cite{metz2020tasks} construct a large and diverse dataset of machine learning optimization tasks, and use it to train optimizers which generalize well.

One drawback of some of the above methods is that the MO has to process extremely high-dimensional data \citep{li2017learning}, rendering its model complexity, training time and deployment unsustainable. 
Some works tried to circumvent this issue by constraining coordinate-wise MOs, but had difficulty in generalizing to new tasks \citep{andrychowicz2016learning}.


\section{Method}
\label{method}
Optimal worst case subspace directions have already been provided by \citep{nemirovski1982orth} for smooth and convex functions.
However, for a given (non worst) case and subject to the limitations of numerical optimization, there might exist better schemes for defining the subspace.
Combining the ideas of meta learning and subspace optimization, we present in this section the proposed Meta Subspace Optimization (MSO) framework. 
We first begin by a general formulation of the framework and then constrain it for better computational efficiency. 

\subsection{A General Framework}
\label{subsec: A general framework}

At each iteration $k \ge 0$, let $S_k = [v_1, v_2, ..., v_L]$ be a list of $L$ vectors spanning a subspace of $\mathbb{R}^n$, where 
$S_{k}^{i}=v_{i}\in \mathbb{R}^{n}$ denotes the $i^{\text{th}}$ element of the list, and $L$ is defined as the memory constraint.
The list $S_k$ comprises first order information, such as current and previous iterates, steps, and gradients.
Let $\Omega_{k}$ denote information obtained from the subspace optimizer at iteration $k-1$, such as subspace iterates \textcolor{black}{$\alpha_k$}, gradients and any by-product of the subspace optimization.
Finally, we define a MO $\pi(S_{k},\Omega_{k})$, which maps $S_k$ and $\Omega_{k}$ to a sequence of $d \le L$ vectors, spanning another subspace of $\mathbb{R}^n$.

In this framework, the subspace minimization problem (\ref{eq:sesop_minimization}) is defined and performed according to the subspace matrix $P_{k}$ obtained from the previous iteration such that \mbox{$P_k = \pi\left(S_{k-1},\Omega_{k-1}\right) \in \mathbb{R}^{n \times d}$}.
Thus, the MO is used to determine at each step $k$, based on $L$ available first order directions \textcolor{black}{and additional information}, a smaller set of directions of size $d$ which is expected to be most beneficial in minimizing the objective in the next step.
Fig. \ref{fig:diagram} shows a diagram of the proposed framework.


\subsection{Meta Subspace Optimization}

Though appealing due to its generality, the MO previously suggested involves the processing of a potentially large number of very high-dimensional input and output vectors, whose dimensions scale with the problem's dimension $n$, \textcolor{black}{which is infeasible for modern problems}.
This issue is reminiscent of other meta-optimization works where some solutions suffer from huge computational cost \citep{li2017learning}, while others employ ad hoc restrictions on the optimization \citep{andrychowicz2016learning}.
Contrary to existing works, in order to make our general framework more efficient for both training and inference, we propose to develop a method which scales with the subspace optimization dimensions $d$, becoming invariant to the size of the original optimization problem $n$.


We first restrict the first order information available to the MO to be only $L$ previously taken steps, such that \mbox{$S_{k}=[p_{t_{L}},\dots,p_{t_{1}}]$}, and we add the \mbox{\emph{mandatory}} current gradient $\nabla f(x_{k})$\citep{nemirovski1982orth} and the current step $p_{k}$ \citep{hestenes1952methods} directly to the subspace matrix $P_{k}$ which we redefine similarly to \citep{conn1996iterated,Narkiss-2005} as follows
\begin{equation}
P_{k}=\big[\pi\big(S_{k-1},\Omega_{k-1}\big),p_{k},\nabla f(x_{k})\big], 
\end{equation}
where \mbox{$t_{L} < ... < t_1 < k$} and $\pi: \mathbb{R}^{n \times L} \rightarrow \mathbb{R}^{n \times (d-2)}$.
During the first $L$ iterations, the list $S_k$ gets populated with the first $L$ steps $p_k$, such that \mbox{$S_{L}=[p_{1},\dots,p_{L}]$}.
In order to maintain a fixed level of complexity (i.e. satisfying the memory constraint $L$), once $S_{k}$ becomes fully populated, for every $k>L$, one of its elements must be removed in order to insert the new step $p_{k+1}$.
Therefore, to further simplify our scheme, we restrict the memory constraint to $L=d-1$, reducing the MO to deciding \emph{which direction should be removed} from $S_{k}$ before adding the new step $p_{k+1}$ at the end of each subspace optimization iteration.
Thus the meta-optimizer is \textcolor{black}{re}defined as $\pi: \mathbb{R}^{n \times (d-1)} \rightarrow \mathbb{R}^{n \times (d-2)}$,
making it more computationally efficient and reducing the \emph{effective} MO's output space to a single dimension.
\textcolor{black}{Note that} under this formulation, \textcolor{black}{one can} add the two other ORTH directions \textcolor{black}{(see Section \ref{subsec:Background-subspace-opt})} in order to enjoy optimal convex worst case convergence of $o(\frac{1}{k^{2}})$ as in \citep{nemirovski1982orth,Narkiss-2005}.

Subspace optimization methods \citep{conn1996iterated, Narkiss-2005,seboost} typically implement a FIFO approach where the oldest subspace direction is removed at each iteration, namely
$\pi({S}_{k},\Omega_{k}) = {S}_{k}\setminus \{{S}_{k}^{1}\}.$
This subspace update implicitly assumes that most of the contribution to the next iterates comes from newly taken steps, whereas the importance of more distant steps taken in the past diminishes.
{\color{black}However, better update strategies exist which are able to remove poor subspace directions from the subspace at earlier stages as follows.}

\subsection{Leveraging Optimal Subspace Step Sizes}
According to the subspace optimization paradigm, at each iteration $k$ a new iterate $x_k$ is updated according to the rule (\ref{eq: subspace optimization update rule}), where $\alpha_k \defeq (\alpha_{k}^{1}, \: \ldots, \: \alpha_{k}^{d}\,)$ is a vector whose elements (henceforth referred to as step sizes) weight each subspace direction. 
Intuitively, 
the vector $\alpha_k$ can be viewed as a multi-dimensional generalization of the standard one-dimensional line search gradient descent.

Therefore, the step sizes obtained from the subspace optimization convey information about the quality of their associated directions in minimizing the objective, where a subspace direction which contributes more to the update of the next iterate will potentially have a larger associated \emph{absolute} step size.
Thus once again we redefine the meta-optimizer as

\begin{equation}\label{eq: general decision rule MSO}
 \pi({S}_{k},\Omega_{k}) := \pi({S}_{k},\alpha_{k}):={S}_{k}\setminus \{{S}_{k}^{
\pi(\alpha_{k})}\},
\end{equation}
where $\pi(\alpha_{k})$ denotes the \emph{effective} meta-subspace optimizer.

This formulation allows very low-dimensional computations since the MO's input dimension scales with the subspace size $d$.
Then, as the optimizer itself, the MO remains invariant to the problem's dimension $n$.
This last assumption suggests an immediate and simple rule-based (RB) approach for updating the subspace by removing at each iteration $k$ the subspace direction $S_k^i$  whose step size $\alpha_k^i$ obtained in the last subspace optimization had the smallest absolute value, i.e.
\begin{equation}\label{eq: decision rule - rule based}
    \pi({S}_{k},\alpha_{k}) = {S}_{k}\setminus \{{S}_{k}^{\argmin |\alpha_{k}|}\}.
\end{equation}
The proposed rule based approach is both simple and efficient and allows improvement over SOTA subspace optimization methods as shown in Fig. \ref{fig:rosenbrock} (left). This approach consistently outperforms the baselines as demonstrated in the Experiments section.

\begin{wrapfigure}{R}{0.52\textwidth}
\vspace{-25pt}
\begin{minipage}{0.52\textwidth}
\begin{algorithm}[H]
\caption{Meta Subspace Optimization (MSO)}\label{algo1}
\begin{algorithmic}
\STATE{\bfseries Input:} $f:\mathbb{R}^{N} \rightarrow \mathbb{R}$, initial point $x_{0}$, meta subspace \\ optimizer \\ $\pi_{\theta}:\mathbb{R}^{(d-1)\times h}\rightarrow \big \{ u \in \mathbb{R}^{d-1} | u^T \mathbf{1} = 1, u_{i}\geq 0 \ \forall i\big \} $ 
\STATE {\bfseries Output:} $x_{final}$
\STATE {\bfseries Initialize:} $P_0 = \{x_0\}$
\FOR{$k=0,1,...,K$}
\IF{$Converged$}
\STATE \textbf{return} $x_{k}$
\ENDIF
\STATE $P_k \leftarrow P_k \cup \{\nabla f(x_k)\}$
\STATE Solve $\alpha_{k}=\argmin_{\alpha}f(x_k+P_{k}\alpha)$ using BFGS
\STATE Update $x_{k+1}=x_{k}+P_{k}\alpha_{k}$
\STATE $P_k \leftarrow P_k \setminus \{\nabla f(x_k)\}$
\IF{$\text{dim}(P_k) \ge d-1$}
\STATE Remove the $a_k^{th}$ direction from $P_k$ where \\ 
$a_k \sim \pi_{\theta}\left(a|\left\{\alpha_{k-h+1}, .., \alpha_{k}\right\}\right)$
\ENDIF
\STATE $p_{k+1} = P_k \alpha_k$
\STATE $P_{k+1} \leftarrow P_k \cup \left\{p_{k+1}\right\}$
\ENDFOR
\STATE \textbf{return} $x_{K}$
\end{algorithmic}
\end{algorithm}
\end{minipage}
\vspace{-15pt}
\end{wrapfigure}
\
\

\subsection{A Reinforcement Learning Approach}
The rule-based approach described above suggests the possibility of finding even better decision rules (or MOs) for selecting which direction to remove from the subspace, based on the last step sizes. 
Moreover, in order to equip the MO with additional information, we allow it to base its decision on the set of vectors $\alpha(k, h)\defeq \left\{\alpha_{k-h+1}, ..., \alpha_k\right\}_{}^{}$ obtained from the last $h$ subspace optimization iterations. 

Based on Eq. (\ref{eq: general decision rule MSO}) we focus our attention on rules of the form:
\begin{equation}\label{eq: general decision rule for RL}
    \pi(S_k, \Omega_{k}) := \pi \left(S_k, \alpha(k, h)\right) := {S}_{k}\setminus \{{S}_{k}^{\pi(\alpha(k, h))}\}
\end{equation}
where $\pi(\alpha(k, h))$ is the index of the susbapce direction to be removed from the subspace.

In order to find superior policies, we resort to reinforcement learning policy search methods. 
Specifically, we model the sequence of decisions - which direction to remove from the subspace - as a Markov Decision Process (MDP). The MDP is defined by the tuple $\left(\mathcal{S}, \mathcal{A}, p\left(s_{j+1}|s_{j}, a_{j}\right), r(s, a)\right)$, where ${\mathcal{S} \subseteq \mathbb{R}^{h\times d}}$ is the state space, corresponding to the $h$ preceding subspace step-size vectors $\left\{\alpha_{k-h+1}, ..., \alpha_k\right\}_{}^{}$,
${\mathcal{A}=\left\{1, ..., d-1\right\}}$ is the action space, corresponding to the $d-1$ subspace directions in $S_{k}=[p_{t_{d-2}},\dots,p_{t_{1}},p_{k}]$, $p\left(s_{j+1}|s_{j}, a_{j}\right)$ is the state transition probability, and $r(s, a)$ is the reward function. 

In this framework, we represent the MO as a stochastic policy \mbox{$\pi_{\theta}(a|s) \defeq \pi\left(a|s;\theta\right)$} parameterized by $\theta$, 
which for a state $s \in \mathcal{S}$, produces a conditional distribution over actions $a \in \mathcal{A}$. 
%
%
%
The Meta Subspace Optimization algorithm is summarized in Algorithm \ref{algo1}. 
At the start of each outer iteration $k$, the subspace $P_k$ includes at most $d-1$ directions, which are all steps taken in the past. 
$P_k$ is then extended with the current gradient $\nabla f(x_k)$ to form a subspace of dimension $d$. 
This subspace is used for solving the subspace optimization (\ref{eq:sesop_minimization}) via the BFGS \citep{NocedWrite} algorithm, yielding the optimal coefficient vector $\alpha_k$. 

The iterate $x_{k+1}$ is then updated according to (\ref{eq: subspace optimization update rule}). 
Next, the gradient $\nabla f(x_k)$ is removed from the subspace.
Then, based on the last $h$ step vectors $\{\alpha_{k-h+1}, ..., \alpha_k\}$, the MO selects which of the $d-1$ steps in the subspace $P_k$ is to be removed.
Finally, the current step $p_{k+1} = x_{k+1} - x_k$ is added to $P_k$ to form $P_{k+1}$ for the next iteration.
Assuming the goal is to learn a MO that minimizes the objective function, the reward $r(s_k,a_k)$ at step $k$ of the optimization is computed based on the relative decrease of the training objective \textcolor{black}{following the subspace optimization} 
and subsequent update of the iterate $x_{k+1}$ such that $r_k \defeq r(s_k, a_k) = (f(x_k)-f(x_{k+1}))/{f(x_k)}.$

In order to train the MO, we define an episode of $T$ steps, during which the MO is used to update the subspace while optimizing an objective function $f(x)$ according to Algorithm \ref{algo1}. 
During training, we use the \mbox{REINFORCE} policy search method \citep{williams1992simple,sutton2018reinforcement} and update the MO at the end of each episode as follows:
$    \theta_{t+1} = \theta_t + \eta
    \frac{1}{m}\sum_{k=1}^{m}\sum_{t=1}^{T} \nabla_{\theta} \log \pi_{\theta}(a_t|s_t) 
    \left(R_t-b(s_t)\right),$
where $\eta$ is the learning rate, $m$ is the batch size, i.e. the number of sampled trajectories used to estimate the policy gradient, ${R_t = \sum_{t'=t}^{T}\gamma^{t'-t}r_{t'}}$ is the return at time step $t$, where $\gamma \in [0,1]$ is the discount rate, and $b(s_t)$ is an estimator of the expected return in state $s_t$ \citep{sutton2018reinforcement}. 
As in many meta-learning tasks, the trained agent is expected to generalize to other tasks, so that the delay incurred by the agent's training is amortized over the inference of new tasks.

\section[title]{Experiments\footnote{Code available at \url{https://github.com/yoniLc/Meta-Subspace-Optimization.git}}}
\paragraph{Experimental Setup}
\label{experiments}
To demonstrate the efficacy of our approach, we train meta subspace optimizers on various non-convex deterministic and stochastic classes of objective functions, namely, the Rosenbrock objective \citep{rosenbrock1960automatic}, robust linear regression \citep{li2017learning}, and neural network classification on the MNIST dataset \citep{lecun-mnisthandwrittendigit-2010}.
We then apply the trained optimizers on new objective functions from the respective classes, and compare their convergence performance to several SOTA subspace optimizers.


\begin{figure}[t!]
 \centering
\begin{subfigure} 
{\includegraphics[trim={0 6 0 6}, width=0.3\textwidth]{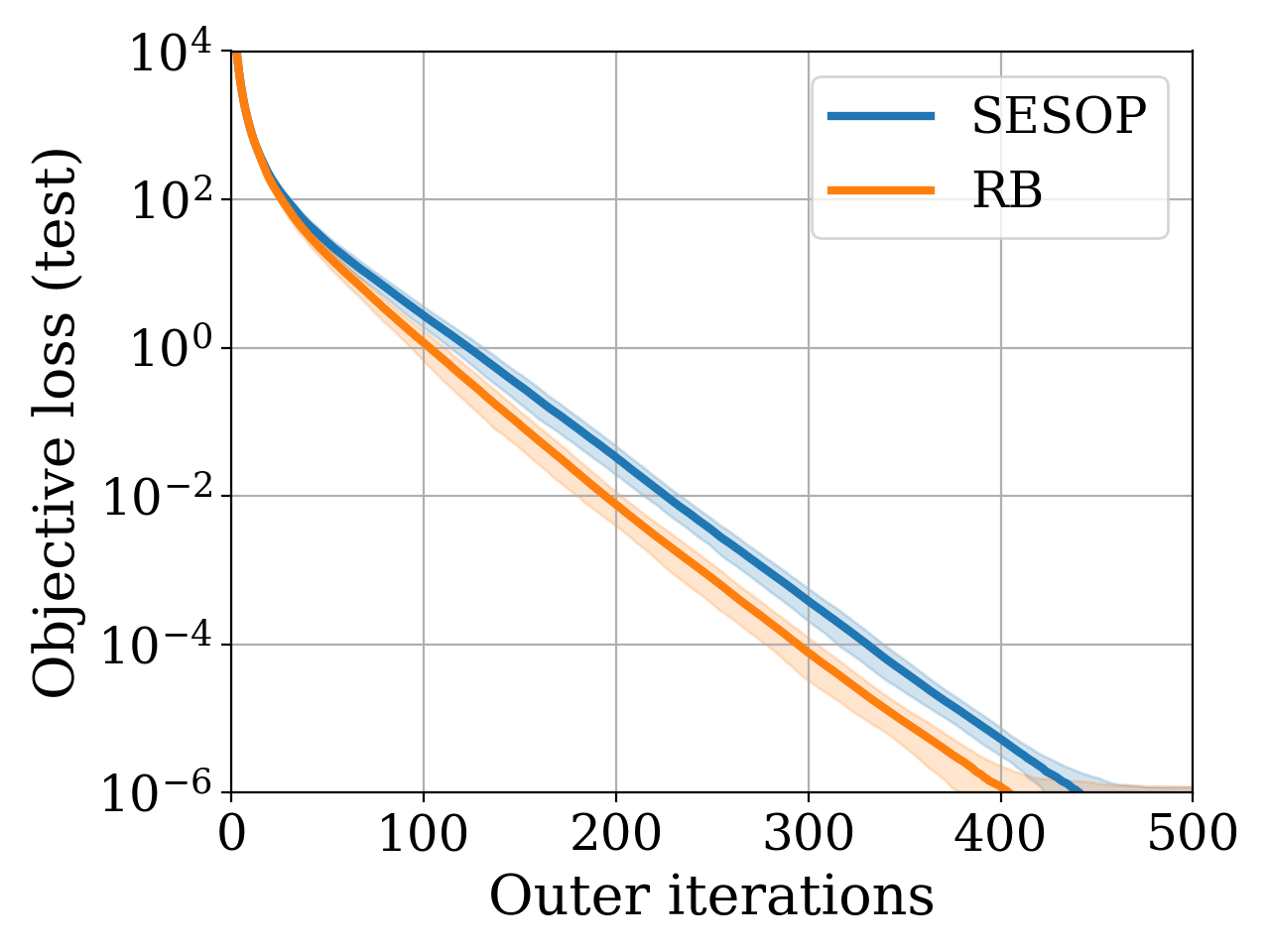}}\hfill
{\includegraphics[trim={0 6 0 6},clip,width=0.3\textwidth]{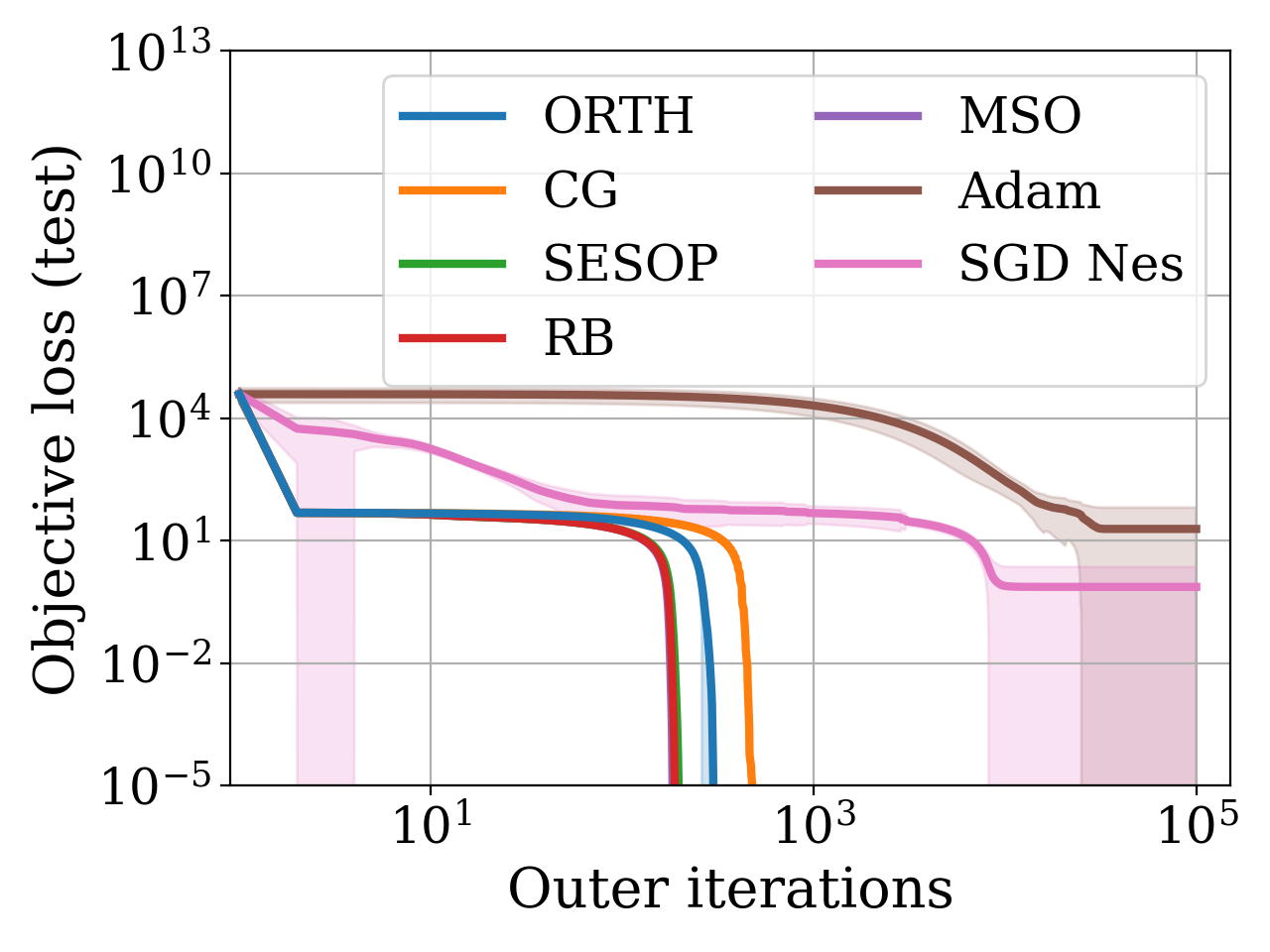}}\hfill
{\includegraphics[trim={13 10 0 15},clip,width=0.3\textwidth]{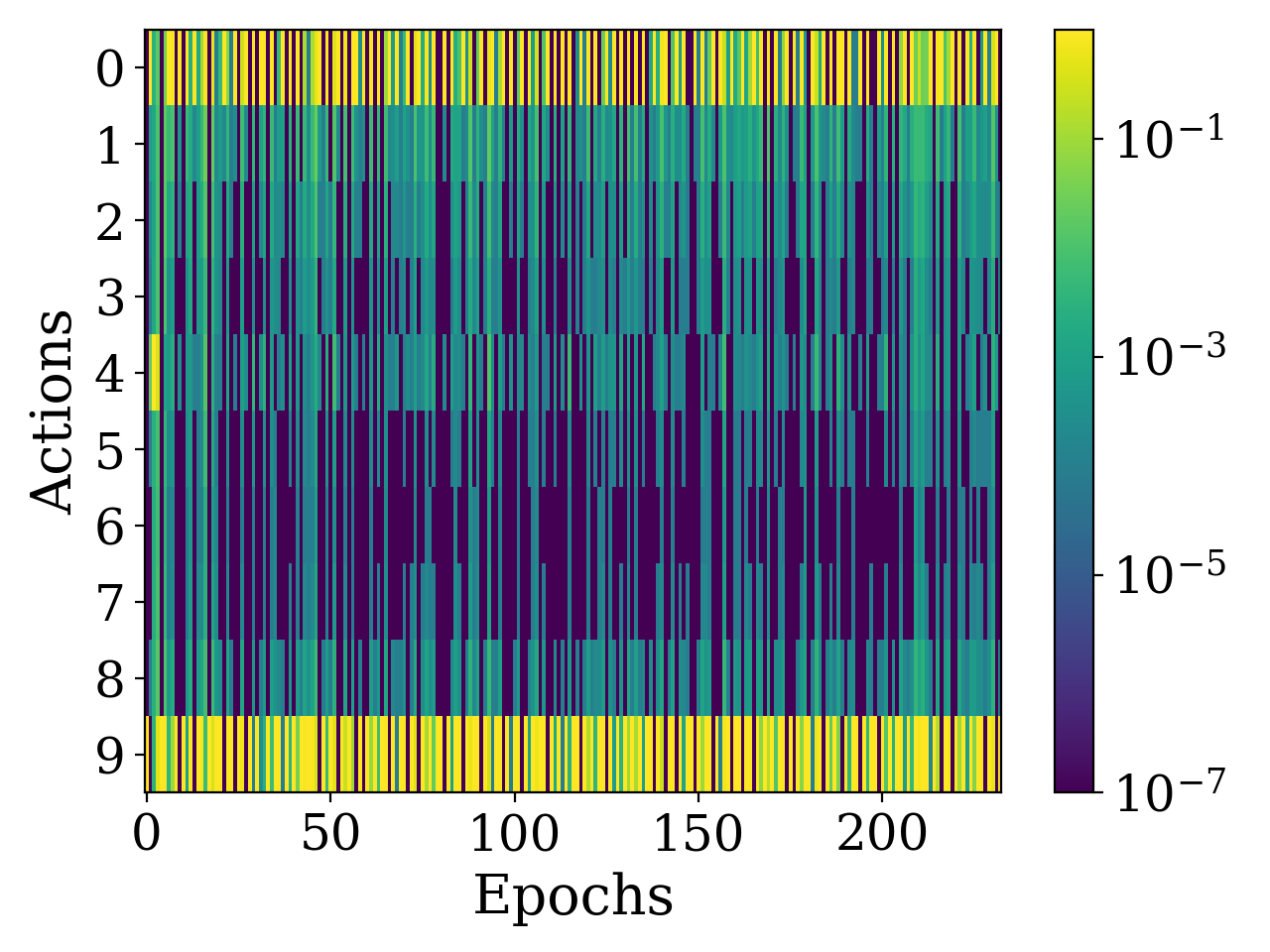}}
\end{subfigure}

\caption{Performance of the baseline subspace optimizer (SESOP) \citep{Narkiss-2005} and the proposed rule based MO (RB) on 100 dimensions quadratic function with $d=10$. Results are averaged over 50 initialization seeds (left). Performance comparison between existing subspace optimizers, the proposed rule-based and the proposed trained agent for the Rosenbrock objective (middle).
MO’s decision trajectory throughout the optimization (right).}
\label{fig:rosenbrock}
\vskip -20pt
\end{figure}

The trained MO is composed of a two-layer fully connected network with 128 neurons in each hidden layer and a tanh nonlinearity, followed by an output layer of dimension $d-1$ and a softmax operation.
The MO is trained using the ADAM optimizer \citep{kingma2014adam} with a learning rate of 5e-3 obtained by grid search. 
The discount rate $\gamma$ is set to 1 to allow consideration of policies whose decisions regarding the subspace have significant impact on future minimization of the objective. Finally, the estimator of the expected return $b(s_t)$ is implemented by an exponential moving average.
\textcolor{black}{Experimental settings for all the different objectives are deferred to the Appendix.}

We compare performance with several known subspace optimization baselines, namely, CG \citep{hestenes1952methods}, ORTH \citep{nemirovski1982orth}, and SESOP \citep{Narkiss-2005}. 
We also report the performance of our rule-based method of Eq. \ref{eq: decision rule - rule based} which we denote \mbox{RB}.
\textcolor{black}{The computational cost of the meta-optimizer neural network is negligible due to its extremely small capacity.}
We redirect the reader to existing literature \citep{Narkiss-2005,conn1996iterated} for comparison of the baselines with other methods (e.g. \citep{Nesterov-1983}  \citep{gill2019practical}).



In all experiments, we run the system depicted in Fig. \ref{fig:diagram}, where the function $f(x)$ is the tested objective function in the deterministic case, or the batch dependent objective in the stochastic case, and the MO module is replaced by the appropriate baseline optimizer where applicable. 
The auxiliary information $\Omega_k$ was set according to the chosen optimizer. 
For instance, in the case of the meta subspace optimizer, $\Omega_k$ is implemented as a matrix whose $i^{th}$ column contains the $h$ last step sizes associated with the $i^{th}$ subspace direction.
In all experiments we use $h=5$.
Step sizes related to subspace directions at iterations for which these directions are not yet available, are set to zero. 
The BFGS algorithm is used with a 1e-5 gradient tolerance threshold. 
In all experiments the subspace dimension is set to $d=10$ as a good trade-off between the optimization performance and the computational cost \citep{seboost} for all the relevant methods. 
\textcolor{black}{Finally, for all subspace optimizers we add the two remaining ORTH directions (see Section \ref{subsec:Background-subspace-opt}) to the subspace at each subspace optimization iteration.}
Time complexity experiment is presented in Figure \ref{fig: analysis_mnist} (right).



\paragraph{Rosenbrock Objective}
In our first experiment, we consider the non-convex Rosenbrock function \citep{rosenbrock1960automatic}: 
$f(x) \defeq \sum_{i=1}^{n-1}\left | b(x_{i+1}-x_i^2)^2 + (a-x_i)^2 \right |$,
where $a, b \in \mathbb{R}$ are constants which modulate the shape of the function and $x_i \in \mathbb{R}$.

We evaluate the MO on 100 instances of this objective, generated using the same procedure as for the training set, and compare its convergence to the baselines. 
As shown in Fig. \ref{fig:rosenbrock} (middle), the MO outperforms all baselines (RB and MSO outperform SESOP by $9\%$ and $14\%$ respectively). 

\paragraph{Robust Linear Regression}
As in \citep{li2017learning}, we consider a linear regression model with a robust loss function.
Training the model requires optimizing the following objective:
$\min_{w,\:b} f(w, b) \defeq \frac{1}{n}\sum_{i=1}^{n}\frac{(y_{i}-w^{T}x_{i}-b)^{2}}{c+(y_{i}-w^{T}x_{i}-b)^{2}},$
where $w\in \mathbb{R}^{D}$ and $b\in \mathbb{R}$ denote the model's weight vector and bias respectively, $x_{i}\in \mathbb{R}^{D}$ and $y_{i}\in \mathbb{R}$ denote the feature vector and label, respectively, of the $i$-th instance, and $c\in \mathbb{R}$ is a constant that modulates the shape of the loss function. 
In our experiments, we use $c = 1$ and $D = 100$. 
This loss function is not convex in either $w$ or $b$.



\begin{figure}[t]
 \centering
\begin{subfigure}
{\includegraphics[trim={0 6 0 0}, width=.3\textwidth,clip]{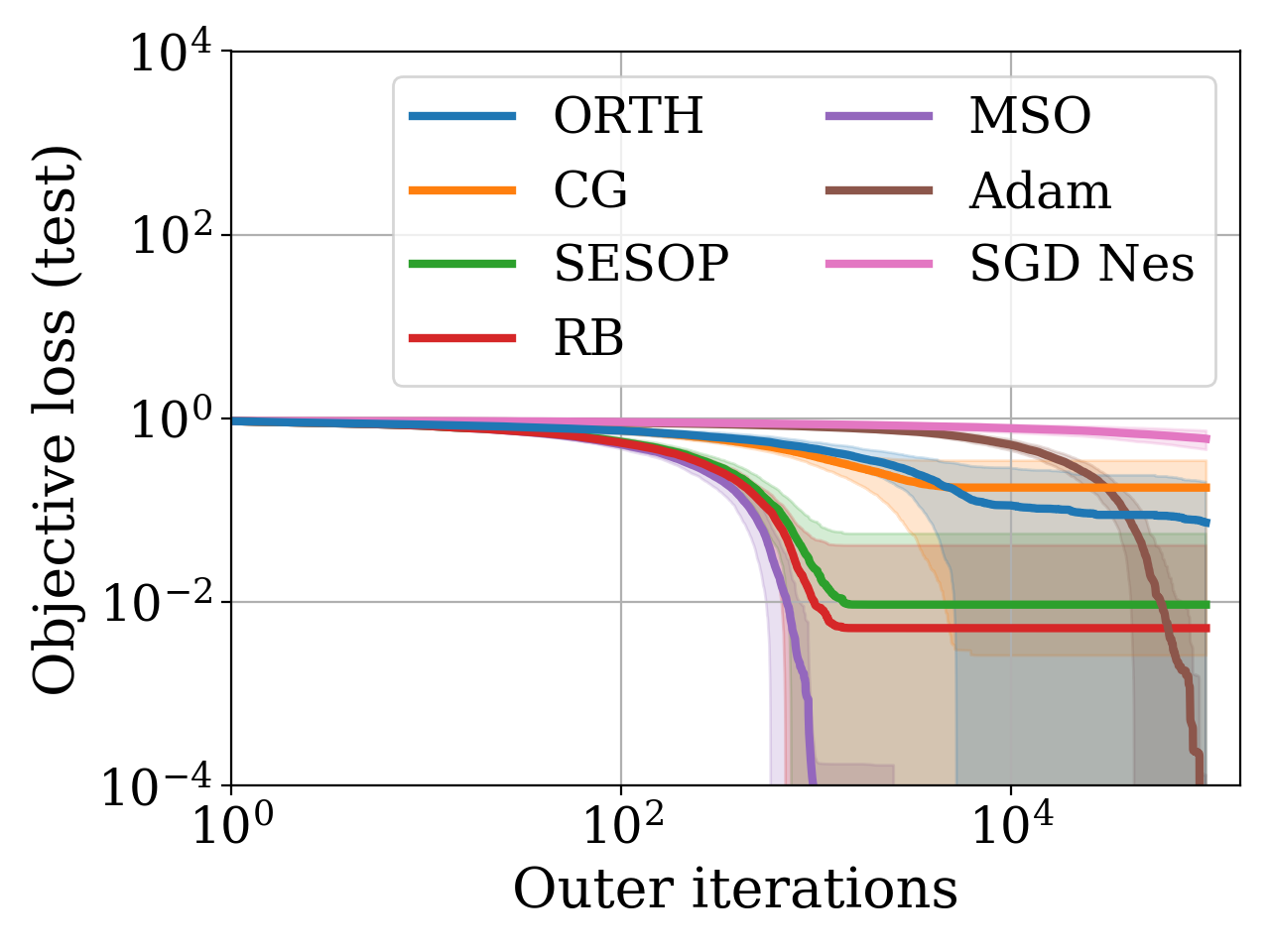}}\hfill
{\includegraphics[trim={13 10 0 15},clip,width=0.3\columnwidth]{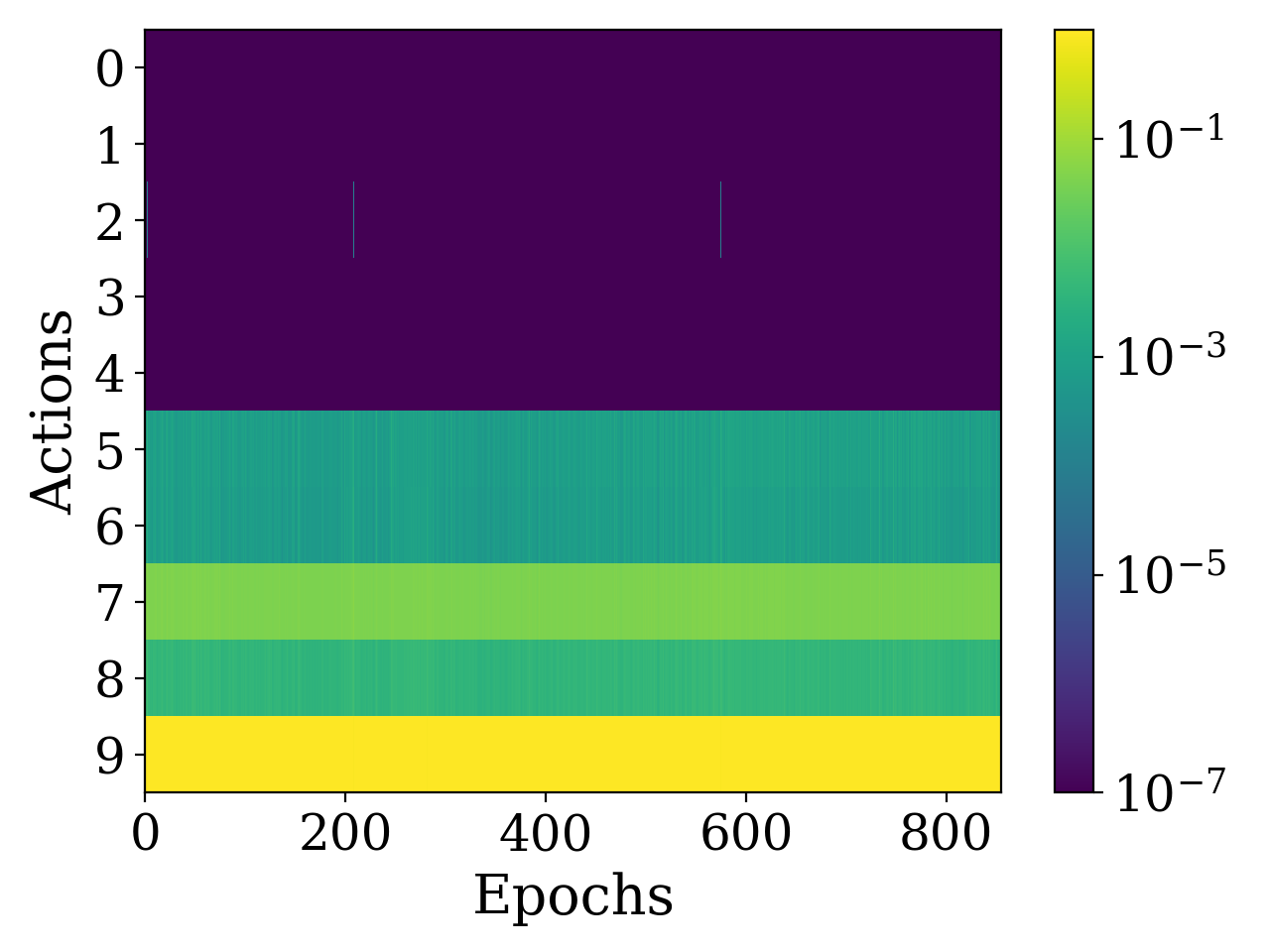}}\hfill
{\includegraphics[trim={0 6 0 0},clip,width=0.3\columnwidth]{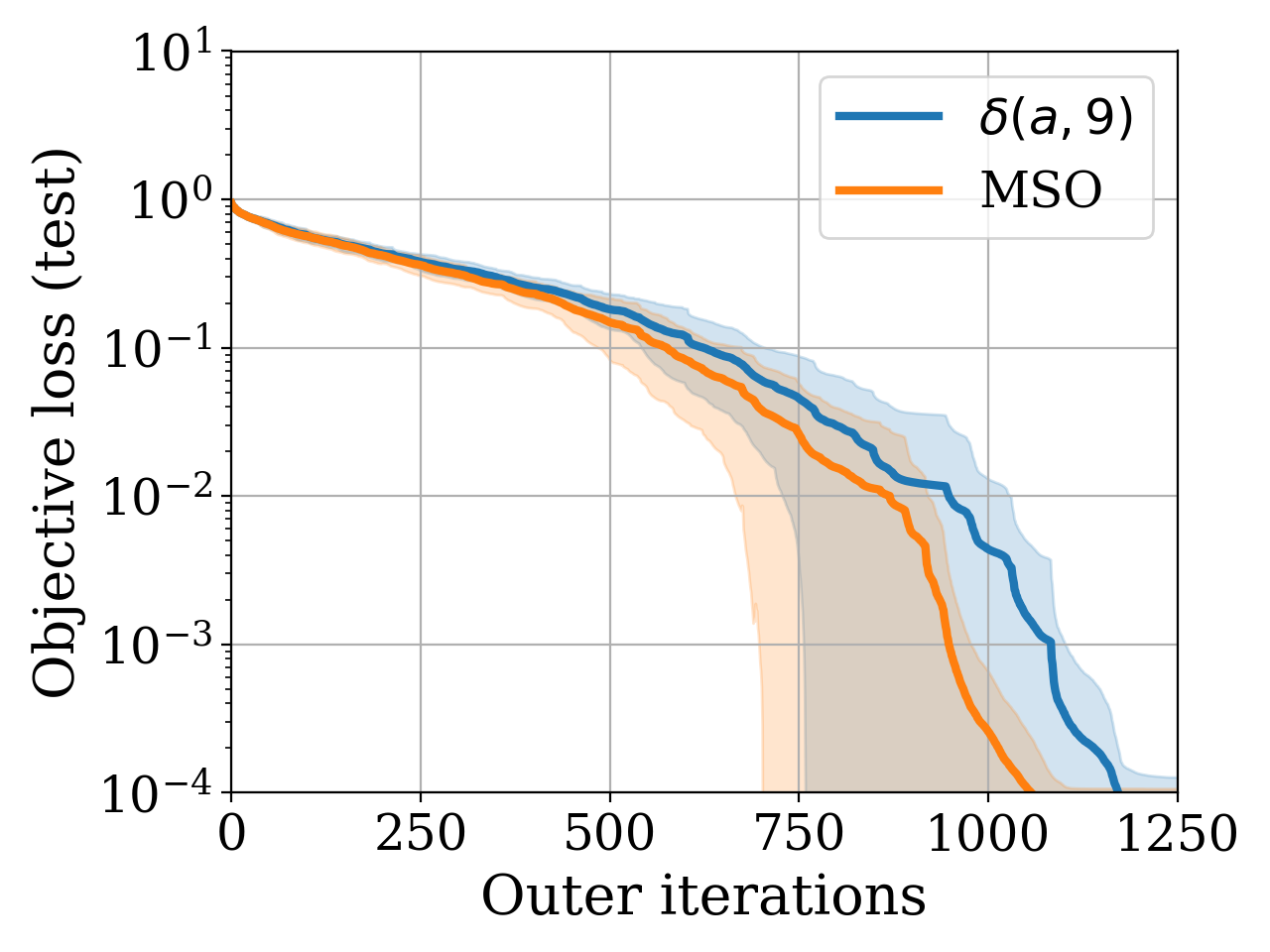}
}
\end{subfigure}
  \caption{Performance of our proposed schemes vs. baselines (left), an example of an MO policy (middle), and the advantage of this policy compared to the deterministic policy agent (right) for the robust linear regression objective.}
\label{fig: robust}
\vskip -10pt
\end{figure}

We evaluate the MO on a test set of 100 random objective functions generated using the same procedure, and compare its convergence to the baselines. 
Fig. \ref{fig: robust} (left) shows for each algorithm the average convergence performance over all test objectives. Some of the baselines suffer from occasional convergence to local minima, which worsen the average convergence curve while our method remains more robust and presents better convergence rate.
Here as well, the meta subspace optimizer outperforms all subspace optimizers throughout the optimization trajectory, even though it was trained only on the first 400 steps. 
Second best is again the \mbox{RB} method which also demonstrates its benefit over traditional subspace optimizers. 

\paragraph{Neural Network Classification}

\begin{figure*}[th!]
\centering
{\includegraphics[trim={0 6 0 0},clip, width=0.30\textwidth]{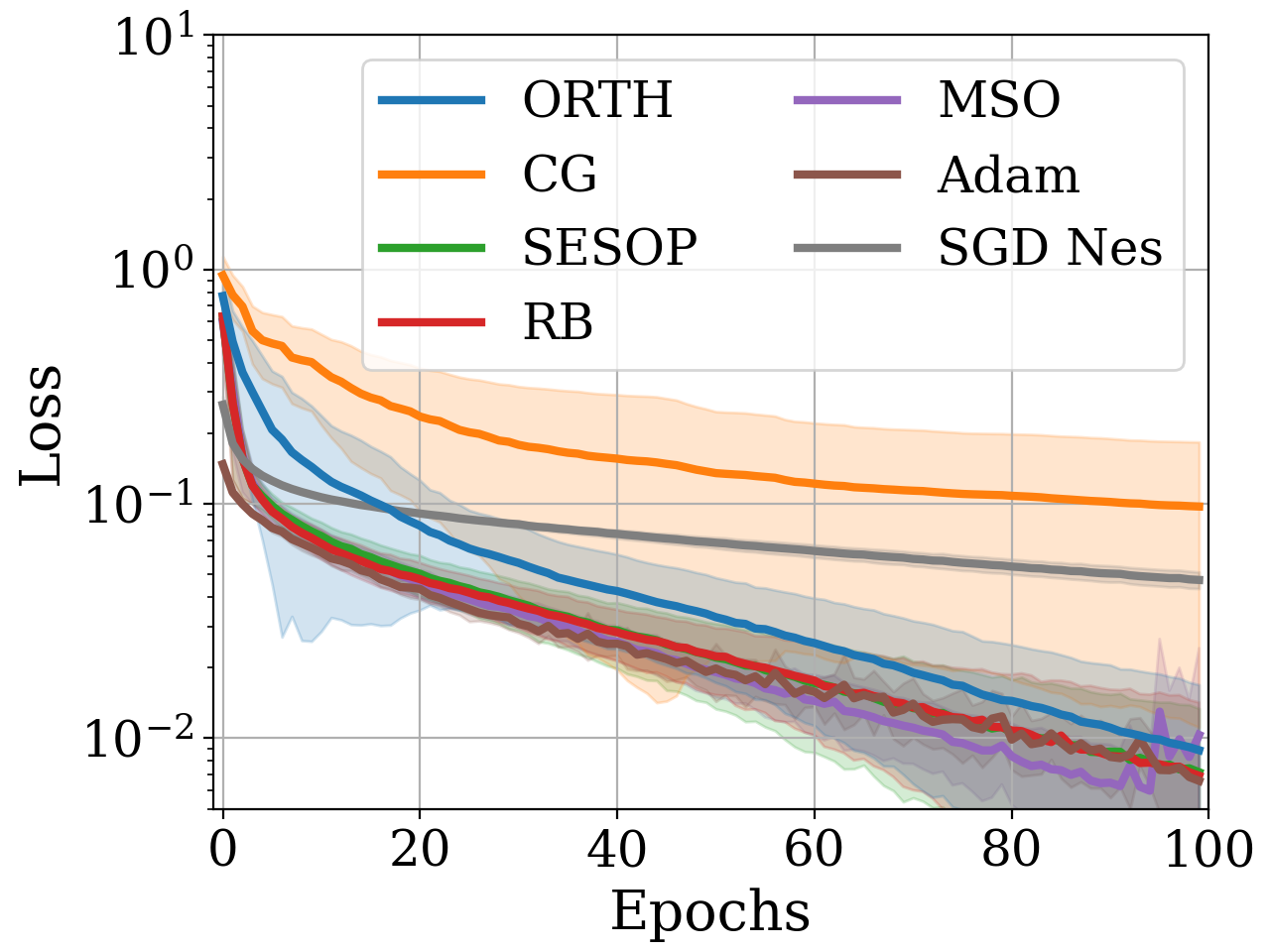}}
{\includegraphics[trim={0 6 0 0},clip, width=0.30\textwidth]{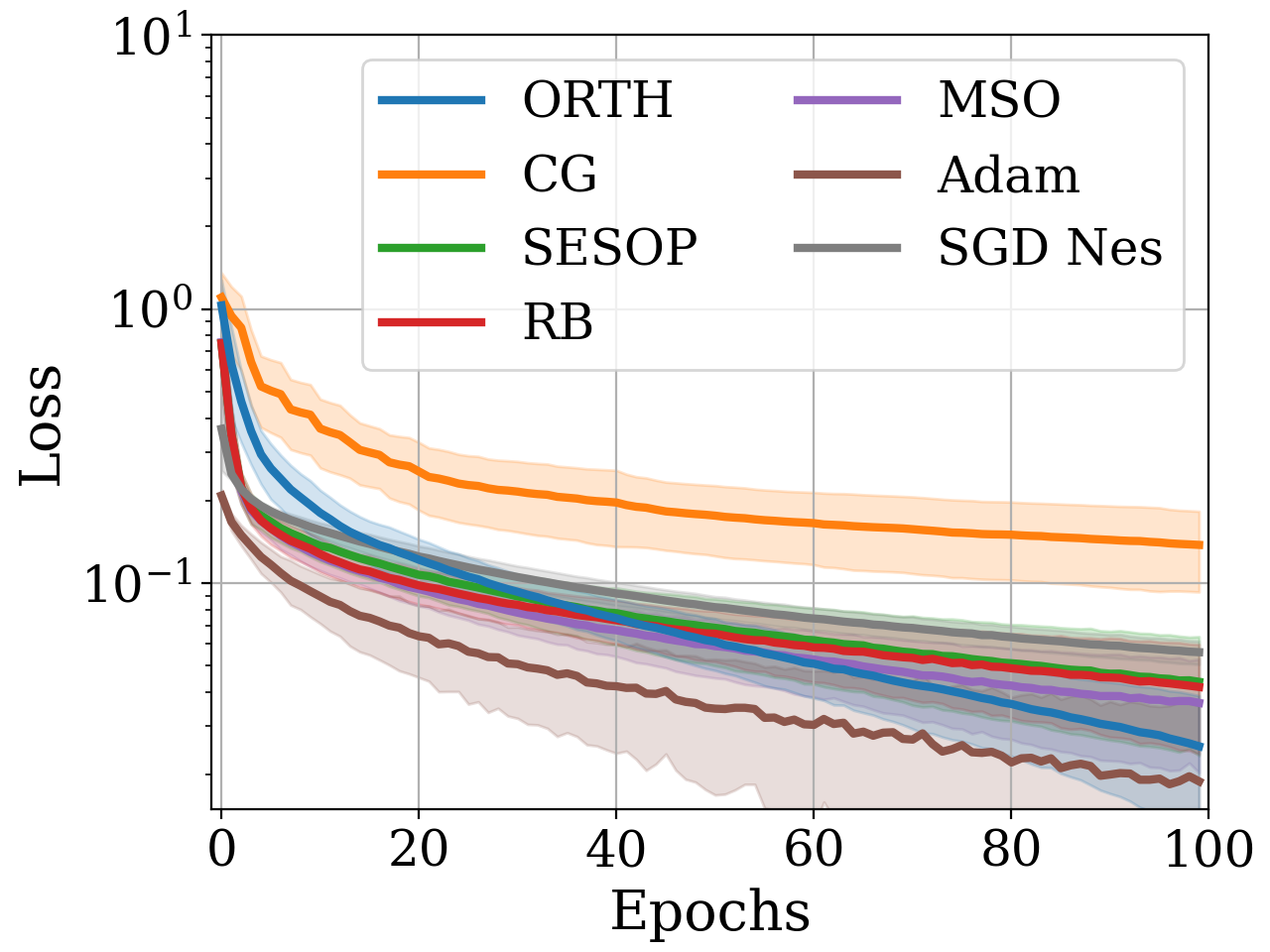}}
{\includegraphics[trim={0 6 0 0},clip, width=0.30\textwidth]{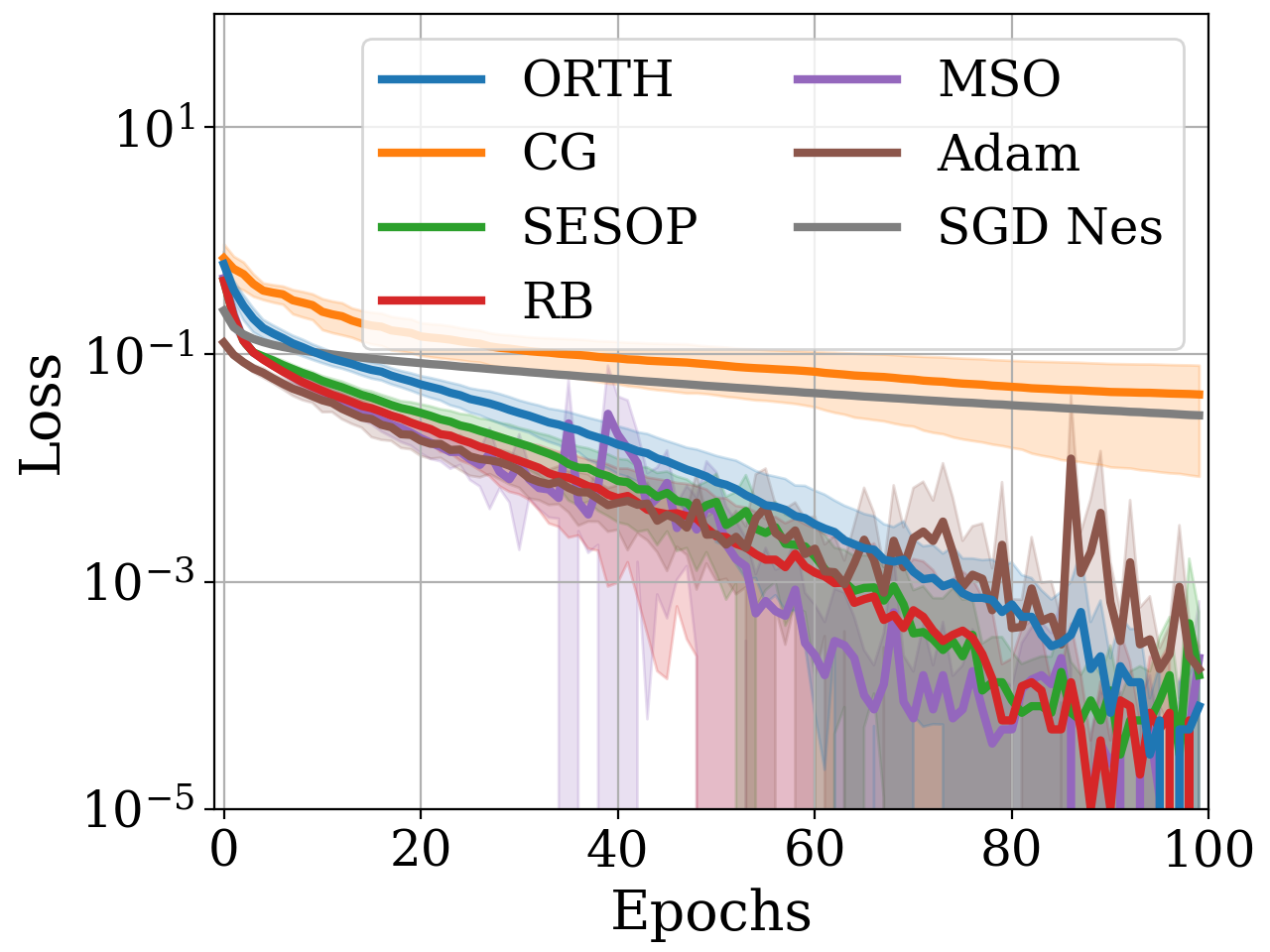}}
\caption{Performance of our proposed schemes vs. baselines when training on half of the MNIST digits (left). Generalization of the MO to the remaining MNIST digits (middle). 
Generalization to a larger network architecture (right). Results are averaged over 10 seeds.}
\label{fig: mnist_exp}
\vskip -20 pt
\end{figure*}


In this experiment, we train a MO for the \emph{stochastic} task of training a small neural network classifier on the MNIST dataset. 
As in \citep{andrychowicz2016learning}, we first train the optimizer to optimize a base network, and then test how well this optimizer generalizes when used to train a neural network with a modified architecture or dataset.

We test two forms of generalization of the MO. 
In a first experiment, we use the trained optimizer to train the base network on the complimentary MNIST digits, which were not used during the MO's training. 
In a second experiment, we apply the MO to train a modified version of the base network with 20 hidden units, instead of 10.
In both cases, shown in Fig. \ref{fig: mnist_exp} (middle and right), the MO generalizes well, even though it is applied in very different conditions than the ones met during its training, and in the second case it even outperforms all hand-crafted baselines.  
\vspace{-10pt}

\section{Analysis}
\begin{figure}[t]
\vskip -15pt
 \centering
{\includegraphics[trim={13 10 0 15},clip,width=0.3\columnwidth]{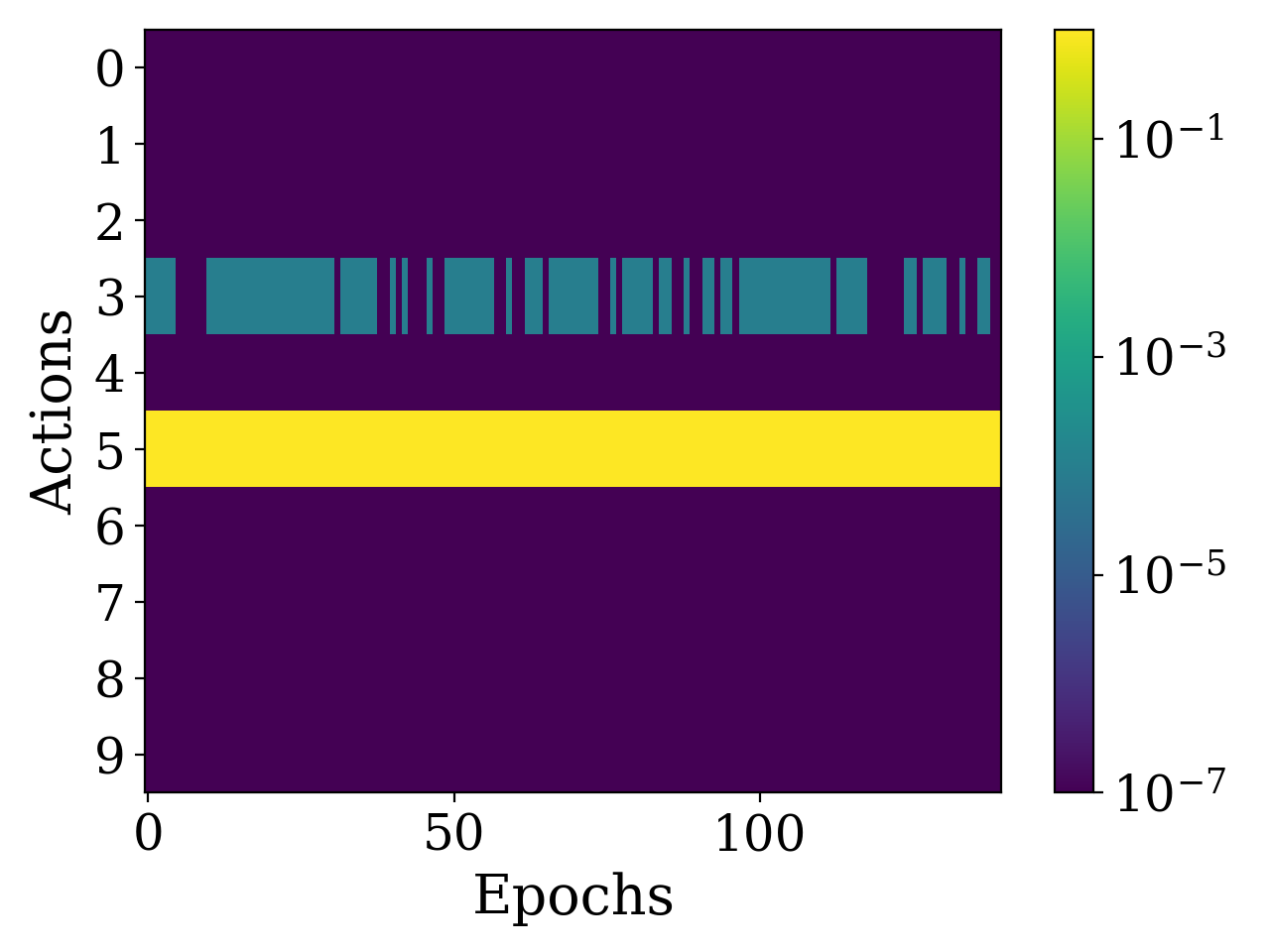}}\hfill
{\includegraphics[trim={0 6 0 6},clip,width=0.3\columnwidth]{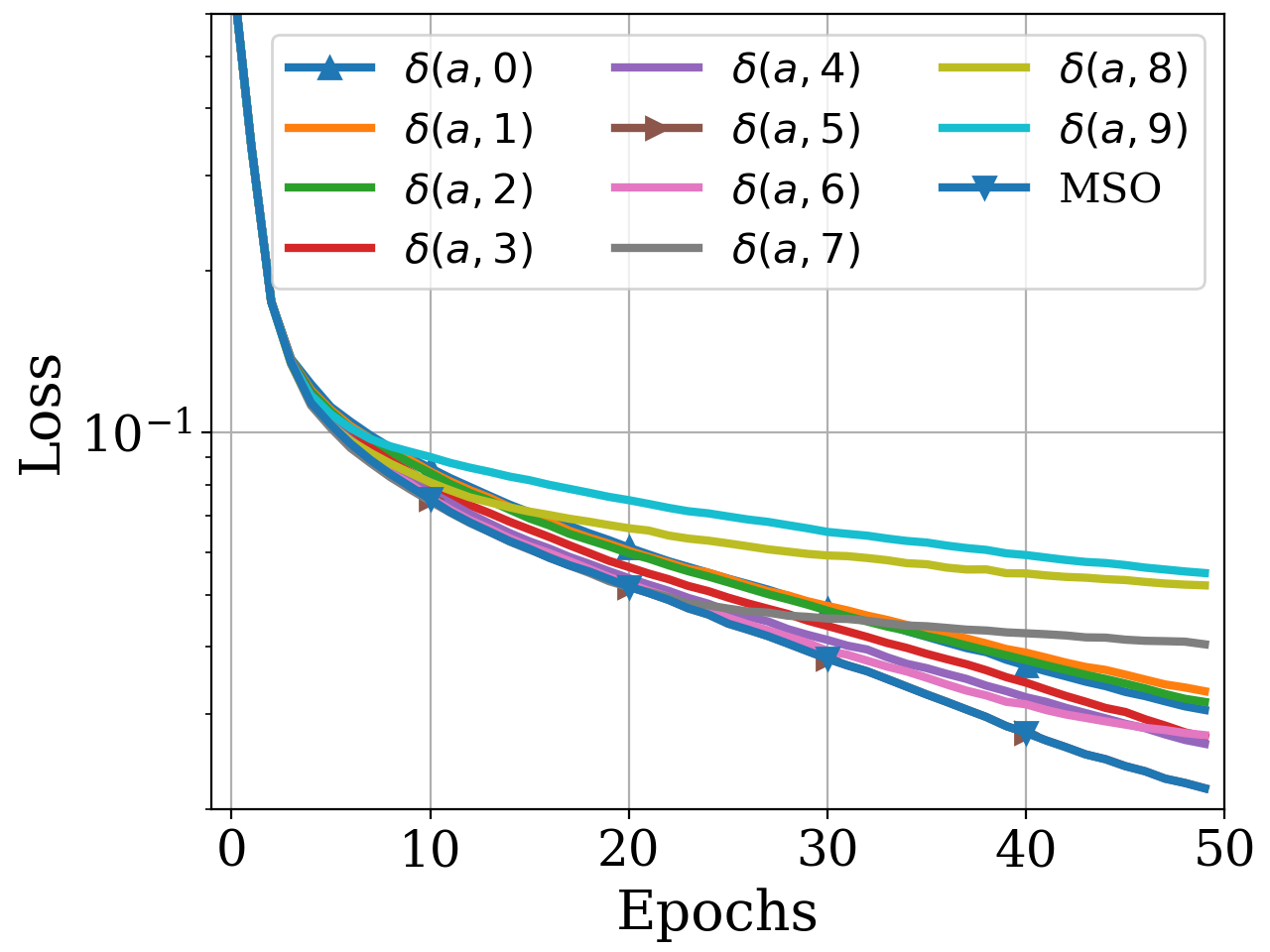}}\hfill
{\includegraphics[trim={0 0 40 20}, width=0.3\columnwidth,clip]{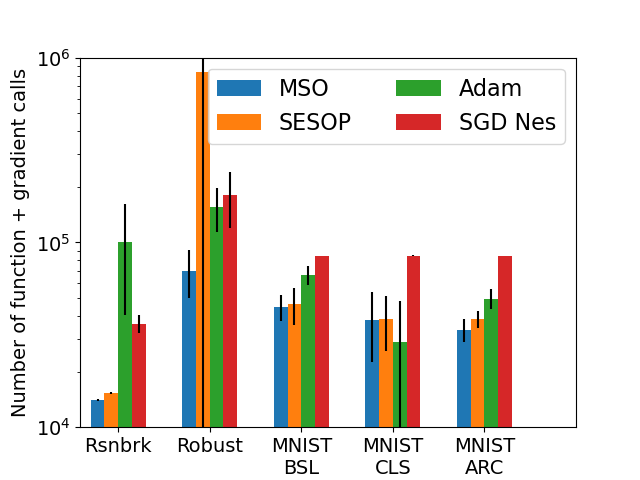}}
\caption{Policy trajectory for MNIST task (left), optimality of the policy w.r.t. all other $\delta$ policies (middle), and comparison of mean and standard variation of the number of function and gradient calls of the methods on the five experiments \textcolor{black}{required for convergence} (right).}
\label{fig: analysis_mnist}
\vskip -20pt
\end{figure}

In this section we analyze the policies learned by the MO, i.e. which subspace direction is least beneficial for minimizing the objective. 
In Fig. \ref{fig:rosenbrock} (right), Fig. \ref{fig: robust} (middle), and Fig. \ref{fig: analysis_mnist} (left) we present examples of the MO's decision trajectory throughout the optimization of the different objectives. 
The $y$ axis represents the action space of the meta optimizer, which in our setting is the set of indices of subspace elements which may be removed from the subspace. For each optimization step, the plots show the probability assigned by the meta optimizer to each of its possible actions. The markedly different policy patterns demonstrate that the MO is able to adapt its policy to each task as we explain next.

In Fig. \ref{fig:rosenbrock} (right) we present the MO's decision trajectory throughout the optimization of the Rosenbrock objective. 
It can be observed that in this setting the MO's policy is highly sensitive to the optimization coefficients $\alpha_k$ associated with the directions of the subspace.
The policy trajectory reveals the following behavior.
Once a new direction $p_k$ is added to the subspace and participates in the optimization for the first time, it is either immediately removed from the subspace, or is allowed to stay in the subspace for subsequent optimizations for as long as it is useful or until it becomes oldest.
This is explained by the large action probability mass located at indices 0 and 9, corresponding to the removal of the oldest and newest directions, respectively, and from the smaller non-zero probabilities for removing the intermediate indices. 

In Fig. \ref{fig: robust} (middle) we present the MO's decision trajectory throughout the optimization of the robust linear regression objective. 
In this case, it can be observed that the MO removes with high probability the \mbox{\emph{most recent}} direction inserted into the subspace matrix, and keeps all remaining subspace directions intact, letting them participate in the next subspace optimization.
Interestingly, this policy is diametrically opposed to existing approaches which axiomatically remove the \emph{oldest} direction from the subspace.
Also, with some (lower) probability, the MO may decide to remove a more distant subspace direction instead, thus \emph{stochastically} allowing the new direction $p_k$ to remain longer in the subspace.

To further asses the impact of this policy, we compare its performance to a \emph{deterministic} policy given by \mbox{$\pi(a|s)=\delta(a,9)$}, where $\delta(n,m)$ denotes the Kronecker delta function.
Such a singular policy always removes the most recent direction from the subspace.
Fig. \ref{fig: robust} (right) shows the average performance of the aforementioned MO compared to a $\delta(a,9)$ policy, demonstrating the benefit of the obtained \emph{stochastic} policy.  
This finding is reminiscent of the improvement obtained with anchor (fixed) points set into the subspace matrix \citep{seboost}.


Fig. \ref{fig: analysis_mnist} (left) presents an example of a MO's decision trajectory applied to the Neural Network (MNIST) classification objective. 
The small probabilities for index 3 imply that this policy is in fact a $\delta(a,5)$ singular policy.
%
%
%
When this policy is applied, a new step $p_k$ which joins the subspace at iteration $k$, remains there for 5 iterations, and is then removed by the MO. 
Moreover, the first 5 steps of the optimization remain in the subspace throughout the course of optimization. 
Thus, the MO equally divides the subspace directions between a set of very early directions which remain relevant even in later optimization steps, and another set of more recent directions which are relevant for a shorter optimization horizon. 

We show that the obtained MO's policy to equally divide the subspace between a long term part and a short term part is optimal among all possible 2-way contiguous subspace partitions.
Fig.~\ref{fig: analysis_mnist}~(middle) compares the performance of the MO to the performance of all singular policies, namely, $\delta(a,0), ..., \delta(a,9)$. 
As can be seen, singular policies which are based on extreme index values are substantially worse than those based on the middle indices, where the optimal singular policy is $\delta(a,5)$ which coincides with the MO's policy.
Apparently, in this setting, devoting half of the subspace for earlier directions in the optimization, and another half for most recent directions gives the best result.
We believe this policy is to be linked with the importance of the initial point $x_{0}$ in the subspace as in the ORTH framework \citep{nemirovski1982orth}. 

Finally, we provide in Fig. \ref{fig: analysis_mnist} (right) the total number of functions and gradients calls of the proposed scheme compared to the SESOP baseline and most importantly compared with the first order methods.
This analysis is to be preferred to timing in order to remain invariant to the implementation of the frameworks.
We can see our approach is superior to the others, especially to first order methods while not requiring hyper-parameter tuning (e.g. learning rate) of the optimizer.
\vspace{-10pt}

\section{Conclusion}
We propose for the first time a meta optimization framework applied to the subspace optimization paradigm.
We derive from our general framework an efficient reinforcement learning approach where the dimensions of the state and action spaces scale with the subspace dimension such that the agent remains invariant to the size of the optimization problem, enabling efficient training and deployment.
We also provide a rule-based method and a reinforcement learning agent able to outperform existing subspace optimization methods, setting a new SOTA subspace optimization method.
Finally, we analyze the obtained agent in order to understand the optimization dynamics.
We believe meta subspace optimization can become a numerical optimization tool of choice for the development of better and more interpretable optimizers.
Future works should pursue overcoming the limitations of the (meta) subspace optimization paradigm on stochastic problems.

\newpage
\bibliographystyle{apalike}
\bibliography{main}

\begin{thebibliography}{}

\bibitem[Andrychowicz et~al., 2016]{andrychowicz2016learning}
Andrychowicz, M., Denil, M., Colmenarejo, S.~G., Hoffman, M.~W., Pfau, D.,
  Schaul, T., Shillingford, B., and de~Freitas, N. (2016).
\newblock Learning to learn by gradient descent by gradient descent.
\newblock In {\em Proceedings of the 30th International Conference on Neural
  Information Processing Systems}, pages 3988--3996.

\bibitem[Bottou et~al., 2018]{bottou2018optimization}
Bottou, L., Curtis, F.~E., and Nocedal, J. (2018).
\newblock Optimization methods for large-scale machine learning.
\newblock {\em Siam Review}, 60(2):223--311.

\bibitem[Choi et~al., 2019]{choi2019empirical}
Choi, D., Shallue, C.~J., Nado, Z., Lee, J., Maddison, C.~J., and Dahl, G.~E.
  (2019).
\newblock On empirical comparisons of optimizers for deep learning.
\newblock {\em arXiv preprint arXiv:1910.05446}.

\bibitem[Choukroun et~al., 2020]{choukroun2020primal}
Choukroun, Y., Zibulevsky, M., and Kisilev, P. (2020).
\newblock Primal-dual sequential subspace optimization for saddle-point
  problems.
\newblock {\em arXiv preprint arXiv:2008.09149}.

\bibitem[Conn et~al., 1996]{conn1996iterated}
Conn, A., Gould, N., Sartenaer, A., and Toint, P.~L. (1996).
\newblock On iterated-subspace minimization methods for nonlinear optimization.
\newblock {\em Linear and Nonlinear Conjugate Gradient-Related Methods}, pages
  50--78.

\bibitem[Cragg and Levy, 1969]{cragg1969study}
Cragg, E. and Levy, A. (1969).
\newblock Study on a supermemory gradient method for the minimization of
  functions.
\newblock {\em Journal of Optimization Theory and Applications}, 4(3):191--205.

\bibitem[Daniel et~al., 2016]{daniel2016learning}
Daniel, C., Taylor, J., and Nowozin, S. (2016).
\newblock Learning step size controllers for robust neural network training.
\newblock In {\em Proceedings of the AAAI Conference on Artificial
  Intelligence}, volume~30.

\bibitem[Dennis~Jr and Turner, 1987]{dennis1987generalized}
Dennis~Jr, J.~E. and Turner, K. (1987).
\newblock Generalized conjugate directions.
\newblock {\em Linear Algebra and its Applications}, 88:187--209.

\bibitem[Gill et~al., 2019]{gill2019practical}
Gill, P.~E., Murray, W., and Wright, M.~H. (2019).
\newblock {\em Practical optimization}.
\newblock SIAM.

\bibitem[Hestenes and Stiefel, 1952]{hestenes1952methods}
Hestenes, M.~R. and Stiefel, E. (1952).
\newblock {\em Methods of conjugate gradients for solving linear systems},
  volume~49.
\newblock NBS Washington, DC.

\bibitem[Hochreiter and Schmidhuber, 1997]{hochreiter1997long}
Hochreiter, S. and Schmidhuber, J. (1997).
\newblock Long short-term memory.
\newblock {\em Neural computation}, 9(8):1735--1780.

\bibitem[Kingma and Ba, 2014]{kingma2014adam}
Kingma, D.~P. and Ba, J. (2014).
\newblock Adam: A method for stochastic optimization.
\newblock {\em arXiv preprint arXiv:1412.6980}.

\bibitem[LeCun and Cortes, 2010]{lecun-mnisthandwrittendigit-2010}
LeCun, Y. and Cortes, C. (2010).
\newblock {MNIST} handwritten digit database.

\bibitem[Li and Malik, 2017]{li2017learning}
Li, K. and Malik, J. (2017).
\newblock Learning to optimize neural nets.
\newblock {\em arXiv preprint arXiv:1703.00441}.

\bibitem[Lv et~al., 2017]{lv2017learning}
Lv, K., Jiang, S., and Li, J. (2017).
\newblock Learning gradient descent: Better generalization and longer horizons.
\newblock In {\em International Conference on Machine Learning}, pages
  2247--2255. PMLR.

\bibitem[Metz et~al., 2020]{metz2020tasks}
Metz, L., Maheswaranathan, N., Freeman, C.~D., Poole, B., and Sohl-Dickstein,
  J. (2020).
\newblock Tasks, stability, architecture, and compute: Training more effective
  learned optimizers, and using them to train themselves.
\newblock {\em arXiv preprint arXiv:2009.11243}.

\bibitem[Metz et~al., 2019]{metz2019understanding}
Metz, L., Maheswaranathan, N., Nixon, J., Freeman, D., and Sohl-Dickstein, J.
  (2019).
\newblock Understanding and correcting pathologies in the training of learned
  optimizers.
\newblock In {\em International Conference on Machine Learning}, pages
  4556--4565. PMLR.

\bibitem[Miele and Cantrell, 1969]{miele1969study}
Miele, A. and Cantrell, J. (1969).
\newblock Study on a memory gradient method for the minimization of functions.
\newblock {\em Journal of Optimization Theory and Applications}, 3(6):459--470.

\bibitem[Narkiss and Zibulevsky, 2005]{Narkiss-2005}
Narkiss, G. and Zibulevsky, M. (2005).
\newblock Sequential subspace optimization method for large-scale unconstrained
  problems.
\newblock Technical Report {CCIT} 559, Technion -- Israel Institute of
  Technology, Faculty of Electrical Engineering.

\bibitem[Nemirovski, 1982]{nemirovski1982orth}
Nemirovski, A. (1982).
\newblock Orth-method for smooth convex optimization.
\newblock {\em Izvestia AN SSSR, Transl.: Eng. Cybern. Soviet J. Comput. Syst.
  Sci}, 2:937--947.

\bibitem[Nesterov, 1983]{Nesterov-1983}
Nesterov, Y. (1983).
\newblock A method for unconstrained convex minimization problem with the rate
  of convergence o(1/n$^2$) (in {Russian}).
\newblock {\em Doklady AN SSSR (the journal is translated to English as Soviet
  Math. Docl.)}, 269(3):543--547.

\bibitem[Nocedal and Wright, 2006]{NocedWrite}
Nocedal, J. and Wright, S. (2006).
\newblock {\em Numerical optimization}.
\newblock Series in operations research and financial engineering, Springer,
  New York,.

\bibitem[Richardson et~al., 2016]{seboost}
Richardson, E., Herskovitz, R., Ginsburg, B., and Zibulevsky, M. (2016).
\newblock Seboost-boosting stochastic learning using subspace optimization
  techniques.
\newblock In {\em Advances in Neural Information Processing Systems}, pages
  1534--1542.

\bibitem[Rosenbrock, 1960]{rosenbrock1960automatic}
Rosenbrock, H. (1960).
\newblock An automatic method for finding the greatest or least value of a
  function.
\newblock {\em The Computer Journal}, 3(3):175--184.

\bibitem[Sutskever et~al., 2013]{sutskever2013importance}
Sutskever, I., Martens, J., Dahl, G., and Hinton, G. (2013).
\newblock On the importance of initialization and momentum in deep learning.
\newblock In {\em International conference on machine learning}, pages
  1139--1147. PMLR.

\bibitem[Sutton and Barto, 2018]{sutton2018reinforcement}
Sutton, R.~S. and Barto, A.~G. (2018).
\newblock {\em Reinforcement learning: An introduction}.
\newblock MIT press.

\bibitem[Vinyals and Povey, 2012]{vinyals2012krylov}
Vinyals, O. and Povey, D. (2012).
\newblock Krylov subspace descent for deep learning.
\newblock In {\em Artificial Intelligence and Statistics}, pages 1261--1268.

\bibitem[Wichrowska et~al., 2017]{wichrowska2017learned}
Wichrowska, O., Maheswaranathan, N., Hoffman, M.~W., Colmenarejo, S.~G., Denil,
  M., Freitas, N., and Sohl-Dickstein, J. (2017).
\newblock Learned optimizers that scale and generalize.
\newblock In {\em International Conference on Machine Learning}, pages
  3751--3760. PMLR.

\bibitem[Williams, 1992]{williams1992simple}
Williams, R.~J. (1992).
\newblock Simple statistical gradient-following algorithms for connectionist
  reinforcement learning.
\newblock {\em Machine learning}, 8(3-4):229--256.

\bibitem[Xu et~al., 2017]{xu2017reinforcement}
Xu, C., Qin, T., Wang, G., and Liu, T.-Y. (2017).
\newblock Reinforcement learning for learning rate control.
\newblock {\em arXiv preprint arXiv:1705.11159}.

\bibitem[Xu et~al., 2019]{xu2019learning}
Xu, Z., Dai, A.~M., Kemp, J., and Metz, L. (2019).
\newblock Learning an adaptive learning rate schedule.
\newblock {\em arXiv preprint arXiv:1909.09712}.

\bibitem[Zibulevsky and Elad, 2010]{ZibEladSPM}
Zibulevsky, M. and Elad, M. (2010).
\newblock {L1-L2} optimization in signal and image processing.
\newblock {\em Signal Processing Magazine, IEEE}, 27(3):76--88.

\end{thebibliography}

\newpage
\appendix
\section{Experimental Settings}
\subsection{Rosenbrock Objective}
We train a MO for optimizing this objective function. 
The training set of the MO consists of 50 instances of this objective differing in their initial point $x_0$ which is randomly drawn from a standard Gaussian distribution. 
In this experiment, the MO is trained for 10,000 episodes where each episode spans a trajectory of 300 optimization steps.
\subsection{Robust Linear Regression}
We train a MO for optimizing objective functions of this form.
The training set (of the MO) consists of 50 examples of such objective functions, where in each example, we randomly draw a dataset of $n=100$ $(x_i, y_i)$ tuples as follows. 
We draw 25 random samples from each one of four multivariate Gaussians, each of which has a random mean and an identity covariance matrix. 
The labels for those points are generated by projecting them along the same random vector, adding a randomly generated bias and perturbing them with i.i.d. Gaussian noise. 

The MO is trained for 2000 episodes. At the start of each episode, the objective function 
is randomly chosen from the training dataset, and the initial values of $w$ and $b$ 
are randomly drawn from a standard Gaussian distribution. 
Each episode is composed of 400 optimization steps.

\subsection{Neural Network Classification}
The objective function to be optimized is the cross entropy loss of a fully connected network with one hidden layer of 10 units and a ReLU nonlinearity, (the algorithm remains stable even with non differentiable activations). 
In order to test generalization of the meta optimizer across different datasets, we first train the meta optimizer to train the base network on a subset of the MNIST dataset consisting half of the digits.
The value of this objective and its gradient are estimated using random minibatches of 8K examples.
The large batch setting is important with second-order optimization in order to obtain good estimate of the true gradient \citep{bottou2018optimization}.
The meta optimizer is trained for 1000 episodes. At the start of each training episode, the dataset is shuffled and the weights are randomly drawn. 
Each episode is run for 14 epochs.
We compare in Fig. 4 
(left) the full batch training loss, where the meta optimizer outperforms the baseline optimizers on this task as well during a period of 100 epochs, after being trained to optimize for only 14 epochs. 

\subsection{Amount of compute ans and resources}
All experiments were conducted on a single Nvidia GeForce GTX 1080 Ti GPU with agent training time ranging between 1-3 days, depending on the task.
Our efficient meta optimizer does not require larger computational resources than the one required for regular first order optimization.




\end{document}